\documentclass[11pt]{amsart}
\usepackage[left=2.0cm,top=2cm,right=2cm]{geometry}
\usepackage{bbm}
\usepackage{mathrsfs}
\usepackage[shortlabels]{enumitem}
\usepackage[usenames,dvipsnames,svgnames,table]{xcolor}
\usepackage[pagebackref,colorlinks=true,linkcolor=BrickRed,citecolor=OliveGreen,pdfstartview=FitH]{hyperref}
\usepackage{graphicx}
\usepackage{float}
\usepackage[normalem]{ulem}
\usepackage{comment}
\usepackage{enumitem}

\newtheorem{thm}{Theorem}[section]

\newtheorem{lem}[thm]{Lemma}
\newtheorem{prop}[thm]{Proposition}

\theoremstyle{definition}

\theoremstyle{remark}

\newtheorem{remk}[thm]{Remark}

\newcommand{\crit}{{\rm crit}}

\newcommand{\Zcal}{ {\mathcal Z}}

\renewcommand{\phi}{\varphi}

\definecolor{blue(ncs)}{rgb}{0.0, 0.53, 0.74}

\title[Schiffer's Conjecture]{A short note on the Schiffer's conjecture for a class of centrally symmetric convex domains in $\mathbb{R}^2$}
\author{Sugata Mondal}
\address{Department of Mathematics and Statistics, University of Reading, UK}
\email{\href{mailto:s.mondal@reading.ac.uk}{s.mondal@reading.ac.uk}}
\date{}

\begin{document}

\maketitle

\begin{abstract}
Let $\Omega$ be a bounded centrally symmetric domain in $\mathbb{R}^2$ with analytic boundary $\partial \Omega$ and center $c$.
Let $\tau = \tau(\Omega)$ be the number of points $p$ on $\partial \Omega$ such that the normal line to $\partial \Omega$ at $p$ passes through $c$. 
We show that if $\tau < 8$ then $\Omega$ satisfies the Schiffer's conjecture.
%\textcolor{blue}{In particular, any sufficiently $C^2$ small analytic perturbation of an ellipse satisfies the conjecture.}
 
%
\end{abstract}

\section{Introduction} 
Let $\Omega$ be a simply-connected bounded domain in the plane with Lipschitz boundary. Assume that there exists a $C^2$ function $u: \Omega \to \mathbb{R}$ that satisfies 
\begin{equation}  \label{eqn:Schiffer}
\Delta u~ =~ \mu \cdot u, ~  \left. \partial_\nu u \right|_{\partial \Omega} \equiv~ 0 \mbox{ \  and \   } \left. u \right|_{\partial \Omega} \equiv~ \mbox{\ 1 \ },
\end{equation}
for some $\mu \neq 0$, where $\partial_\nu$ denotes the unit outward normal vector field along $\partial \Omega$.
The Schiffer's conjecture says that $\Omega$ is a disc \cite{Y92}.
Interest in this conjecture, for bounded connected domains in the plane with analytic boundary, partially, comes from its connection to a well-studied problem in integral geometry known as the Pompieu problem. 
A plane domain $\Omega$ that is not a disc is said to have the Pompieu property if there exists a continuous function $f: \mathbb{R}^2 \to \mathbb{R}$ and a rigid motion $\sigma$ of $\mathbb{R}^2$ such that $$\int_{\sigma(\Omega)} f = 0.$$
 
In \cite{BST} it was proved that if a bounded domain $\Omega$ in the plane does not have the Pompieu property then there exists a solution to \eqref{eqn:Schiffer}.
In \cite{W2} Williams proved a free boundary result
concerning \eqref{eqn:Schiffer}. 
As a consequence of it, in $\mathbb{R}^2$ one obtains that a bounded, simply-connected open set $\Omega \subset \mathbb{R}^2$ with Lipschitz boundary $\partial \Omega $ has the Pompeiu property unless $\partial \Omega$ is analytic.
Hence, to prove that discs are the only bounded simply-connected plane domains with Lipschitz boundary that does not have the Pompieu property, it suffices to prove the Schiffer's conjecture for domains in the plane with analytic boundary. 
Perhaps not entirely unexpectedly, the Schiffer's conjecture, in some form, shows up in the asymptotic properties of
nodal sets of eigenfunctions \cite{TZ}.

For the rest of the paper we will be concerned with simply connected domains in the plane with analytic boundary, unless otherwise stated.

Althgough the literature on the Pompieu problem is fairly rich (see \cite{GF1}, \cite{GF2}, \cite{GF3}, \cite{GF4}, \cite{E1}, \cite{E2}, \cite{E3}, \cite{D}), the Schiffer's conjecture has not been studied that extensively. In \cite{Ber} Bernstein showed that if there are infinitely many solutions to \eqref{eqn:Schiffer}then $\Omega$ must be a disc. In \cite{BY} it was proved that if the $\mu$ in \eqref{eqn:Schiffer} is equal to $\mu_2(\Omega)$ -- the second Neumann eigenvalue of $\Omega$ -- then $\Omega$ is a disc. In \cite{Av} it was shown that if $\Omega$ is convex and if $\mu$ in \eqref{eqn:Schiffer} is at most $\mu_7(\Omega)$ -- the seventh Neumann eigenvalue of $\Omega$ -- then $\Omega$ is a disc.
Recently, Deng \cite{De} has obtained results similar to the last one where he was able to replace the second Neumann eigenvalue of $\Omega$ by some larger eigenvalue.

Other recent developments around the Schiffer's conjecture include \cite{L}, \cite{NSY} and \cite{KM}. 
In this short note we extend an observation of Deng \cite{De} to prove Schiffer's conjecture for a class of centrally symmetric domains. 
The main result of this paper is the following theorem.
\begin{thm}\label{main_1}
Let $\Omega$ be a centrally symmetric domain with analytic boundary $\partial \Omega$ and centre $c$.
Let $\tau = \tau(\Omega)$ be the number of points $p$ on $\partial \Omega$ such that the normal line to $\partial \Omega$ at $p$ passes through $c$. 
We show that if $\tau < 8$ then $\Omega$ satisfies the Schiffer's conjecture. 
\end{thm}
\begin{remk}
A few remarks are in order.
\begin{enumerate}
    \item It can be easily checked that for $\Omega$ an ellipse, the number $\tau =4$. Hence, Theorem \ref{main_1} implies that any ellipse satisfies the Schiffer's conjecture. Because the condition in Theorem \ref{main_1} is `open' it follows that any controlled deformation of any given ellipse (that satisfies the assumptions in Theorem \ref{main_1}) also satisfies the Schiffer's conjecture.

    \item According to \cite{D-L}, for any centrally symmetric convex (c-s-c) domain $\Omega$ in the plane with smooth boundary, and for any $p \in \Omega$, the average number of points on $\partial \Omega$ where the normal line to $\partial \Omega$ passes through $p$, is at most $8$. It is not very difficult to show the existence of c-s-c domains in the plane with $\tau > 8$ (here $\tau$ is as in Theorem \ref{main_1}). It would be interesting to see if one can work with a general point in $\Omega$ rather than the centre and conclude the conjecture for all c-s-c domains from the result of \cite{D-L}. 
\end{enumerate}

\end{remk}

\subsection*{Sketch of the proof of Theorem \ref{main_1}}
%
%In \S2 we provide the details of some results due to Jian Deng \cite{De}. 
Let $\Omega$ be a simply connected domain in the plane with analytic boundary and let $u$ be a solution of \eqref{eqn:Schiffer}.
In \S2.3 we obtain certain integral inequalities that $u$ satisfy. 
When $\Omega$ is centrally symmetric, using these identities, we deduce that the normal derivative along $\partial \Omega$, of a particular {\it rotational} (see \S2.1) derivative of $u$, must vanish at at least eight points (Theorem \ref{thm:deng}).

In \S3 we study local behavior of $u$ on $\partial \Omega$. 
We first obtain a local expression of $u$ at a point on $\partial \Omega$ 
 (Proposition \ref{prop:degree_along_normal}). 
Then, we use this expression to give a geometric description of the points where the normal derivative of a rotational derivative of $u$ may vanish. 
The proof of Theorem \ref{main_1} is deduced from this geometric description.
\subsection*{Acknowldegement}
The author would like to thank Michael Levitin for all the discussions on Schiffer's conjecture, and, for pointing out the very interesting literature around the average number of normals for planar convex domains including \cite{D-L}.
\section{Counting nodal critical points of derivatives of $u$ after Jian Deng}
\subsection{Killing fields}\label{subsection:killingfields}
A Killing field on $\mathbb{R}^2$ is a vector field whose infinitesimal generators are isometries of $\mathbb{R}^2$. 
It is well known that any Killing field on $\mathbb{R}^2$ is either a constant vector field or a rotational vector field.
A constant vector field $L$ is a vector field on $\mathbb{R}^2$ such that there exists $a, b \in \mathbb{R}$ with $a \cdot b \neq 0$ such that $$L = a \cdot \partial_x + b \cdot \partial_y.$$  
Here $x, y$ denote the cartesian coordinates of $\mathbb{R}^2$.
A rotational vector field $R_p$ with centre $p$, on the other hand, is a vector field such that $$R_p = \partial_\theta$$ where $(r, \theta)$ are the polar coordinates on $\mathbb{R}^2$ with centre $p$. 
One property of a Killing field $X$, that we will be using in this paper, is that it commutes with the Laplacian $$\Delta \cdot X = X \cdot \Delta.$$
This property implies that if $v : \mathbb{R}^2 \to \mathbb{R}$ be a $C^3$ function such that $\Delta v = \lambda \cdot v$ on a domain $\Omega$ then $\Delta Xv = \lambda \cdot Xv$ on $\Omega$. 
\subsection{Boundary behavior of solutions of \eqref{eqn:Schiffer}}
Let $\Omega$ be a simply connected domain with analytic boundary $\partial \Omega$.
We denote the unit tangent and the unit outward normal vector fields along $\partial \Omega$ by $\partial_\tau$ and $\partial_\nu$ respectively.

Let $v$ be a Neumann eigenfunction on $\Omega$. 
This means, there exists a $\lambda > 0$ such that $\Delta v = \lambda \cdot v$ and $v$ satisfies Neumann boundary condition $\partial_\nu v = 0$ along $\partial \Omega$.
It is well known that $v$ can be extended analytically to a neighborhood $N(\Omega)$ of $\Omega$ \cite{TZ}.
From here onwards we shall always think of $v$ as function on $N(\Omega)$.
By the analyticity, the equality $\Delta v = \lambda \cdot v$ still holds on all of $N(\Omega)$.

For a function $f: \Omega \to \mathbb{R}$, the set $f^{-1}(0)$ of zeros of $f$ is called the {\it nodal set} of $f$ and we denote it by $\mathcal{Z}(f)$. 
If $f$ is $C^1$, we call a point $p \in \Omega$, a {\it critical point} of $f$ if the gradient of $f$ vanishes at $p$.
Equivalently, $p$ is a critical point of $f$ if all partial derivatives of $f$ at $p$ vanish.
\begin{lem}\label{lem:boundary-crit}
Let $u$ be a solution of \eqref{eqn:Schiffer} on $\Omega$.
Then each point on $\partial \Omega$ is a critical point of $u$.
In particular, $Xu|_{\partial \Omega} \equiv 0$ for any vector field $X$. 
\end{lem}
\begin{proof}
Because $\partial_\tau(p)$ and $\partial_\nu(p)$ are linearly independent, and because $\partial_\nu u|_{\partial \Omega} \equiv 0$, it suffices to show that $\partial_\tau u(p) = 0$ for each $p \in \partial \Omega$.
This follows from the fact that $u \equiv 1$ on $\partial \Omega$.
\end{proof}
Let $s$ denote the arc length parametrization of $\partial \Omega$.
Let the parametric equation of $\partial \Omega$ be given by $z(s) = (x(s), y(s))$.
In the complex notation $z(s) = x(s) + i \cdot y(s)$ such that the derivative $$\frac{dz}{ds} = e^{i \cdot \theta(s)},$$
where $\theta(s)$ is the angle the tangent to $\partial \Omega$ at the point $z(s)$ makes with the $x$ axis.
%The curvature of $\partial \Omega$ at $z(s)$ is given by $\rho(s) = - d\theta/ds$.

\subsubsection{Boundary values of derivatives of $u$}
Most of the results in this subsection are reformulations of some results from \S2 of \cite{De}. We give the details here for the sake of completeness.

The main goal here is to find expressions for $u_{xx}, u_{xy}$ and $u_{yy}$ along $\partial \Omega$, in terms of the angle function $\theta(s)$ defined above.
For this, we first recall that the unit tangent vector to $\partial \Omega$ at $z(s)$ is given by $$\partial_\tau = \cos{\theta(s)} \cdot \frac{\partial}{\partial x} + \sin{\theta(s)} \cdot \frac{\partial}{\partial y}.$$
Also, after possibly reflecting $\Omega$ along the $x$-axis, the unit outward normal vector to $\partial \Omega$ at $z(s)$ is given by $$\partial_\nu = \sin{\theta(s)} \cdot \frac{\partial}{\partial x} - \cos{\theta(s)} \cdot \frac{\partial}{\partial y}.$$ 

Using Lemma \ref{lem:boundary-crit}, we get
\begin{equation}
    0 = \partial_\tau u_x = \cos{\theta(s)} \cdot \frac{\partial^2 u}{\partial x^2} + \sin{\theta(s)} \cdot \frac{\partial^2 u}{\partial x \partial y},
\end{equation}
and, 
\begin{equation}
    0 = \partial_\tau u_y = \cos{\theta(s)} \cdot \frac{\partial^2 u}{\partial x \partial y} + \sin{\theta(s)} \cdot \frac{\partial^2 u}{\partial y^2}.
\end{equation}
Since $$\frac{\partial^2 u}{\partial x^2} + \frac{\partial^2 u}{\partial y^2} = - \mu \cdot u,$$
a straight forward computation shows that 
\begin{equation}\label{uxx}
    u_{xx}|_{\partial \Omega} = \frac{1}{2} \cdot (1 - \cos(2\theta(s))) \cdot \mu \cdot u|_{\partial \Omega}. 
\end{equation}
Similarly, we obtain other identities: $u_{xy}|_{\partial \Omega} = - \frac{1}{2} \cdot \sin(2\theta(s)) \cdot \mu \cdot u|_{\partial \Omega}$, and, $u_{yy}|_{\partial \Omega} = \frac{1}{2} \cdot (1 + \cos(2\theta(s))) \cdot \mu \cdot u|_{\partial \Omega}.$

\subsection{Counting the nodal critical point on the boundary}
Let $X$ be a Killing field on $\mathbb{R}^2$.
Hence, from \S \ref{subsection:killingfields}, either $X = c_1 \cdot \partial_x +  c_2 \cdot \partial_y$, for some constants $c_1, c_2$ (with $c_1 \cdot c_2 \neq 0$), or, $X = R_p$ for some $p \in \mathbb{R}^2$.
By Lemma \ref{lem:boundary-crit}, $\partial \Omega \subset \Zcal(X u)$.
Because $X$ is a Killing field, it follows that $\Delta Xu = \mu \cdot Xu$. Hence, by Green's formula \cite{Cha}
\begin{equation}
    0 = \int_{\partial \Omega} X u \cdot \frac{\partial u_{xx}}{\partial n} ds = \int_{\partial \Omega} \frac{\partial X u}{\partial n} \cdot u_{xx} ds.
\end{equation}
Now, by \eqref{uxx}, we get our first equality
\begin{equation*}
    \int_{\partial \Omega} \frac{\partial X u}{\partial n} \cdot \frac{1}{2} \cdot (1 - \cos(2\theta(s))) \cdot \mu \cdot u|_{\partial \Omega}  ds = 0.
\end{equation*}
Since $u|_{\partial \Omega}$ is a constant, we get
\begin{equation}
    \int_{\partial \Omega} \frac{\partial X u}{\partial n} \cdot \frac{1}{2} \cdot (1 - \cos(2\theta(s))) ds = 0.
\end{equation}
Similarly, working with $u_{xy}$ and $u_{yy}$ we get the identities
\begin{equation}
    \int_{\partial \Omega} \frac{\partial X u}{\partial n} \cdot \sin(2\theta(s)) ds = 0 = \int_{\partial \Omega} \frac{\partial X u}{\partial n} \cdot (1 + \cos(2\theta(s))) ds.
\end{equation}
Combining these equations we get the identities 
\begin{equation}\label{vanishing_final}
    \int_{\partial \Omega} \frac{\partial X u}{\partial n} ~ ds = 0, ~~ \int_{\partial \Omega} \frac{\partial X u}{\partial n} \cdot \sin(2\theta(s)) ~ ds = 0, ~~ \textrm{and} ~~ \int_{\partial \Omega} \frac{\partial X u}{\partial n} \cdot \cos(2\theta(s)) ~ ds = 0.
\end{equation}
Now we use these identities to obtain a lower bound on the number of nodal critical points of $R_c u$ on $\partial \Omega$.
Observe that, because $R_c u \equiv 0$ on $\partial \Omega$, the tangential derivative $\partial_\tau R_c u(p) = 0$ for every $p \in \partial \Omega$.
Hence, every point $p \in \partial \Omega$ where $\partial_\nu R_c u(p) = 0$, is a nodal critical point of $R_c$.
Because $R_c u$ satisfies $\Delta u = \mu \cdot u$ in $N(\Omega)$, the nodal critical points of $R_c$ are isolated by \cite{Cheng}.
In particular, there are only finitely many points on $\partial \Omega$ where $\partial_\nu R_c$ vanishes.
\section{Centrally symmetric convex domains}
Let $\Omega$ be a centrally symmetric convex domain.

\begin{thm}[Deng]\label{thm:deng}
Let $\Omega$ be a centrally symmetric domain with centre $c$. Let $u$ be a solution of \eqref{eqn:Schiffer}. Then there are at least eight points on $\partial \Omega$ where $\partial_\nu R_c$ vanishes.
\end{thm}
\begin{proof}
We first claim that $u$ is centrally symmetric with respect the central symmetry $\iota$ of $\Omega$.
To see this, we consider the function $v = (u + \iota(u))/2$.
Clearly $v$ satisfies \eqref{eqn:Schiffer}.
If $v$ is not a multiple of $u$ then we may consider a linear combination $w$, of $u$ and $v$, such that $w$ satisfies \eqref{eqn:Schiffer} and $w|_{\partial \Omega} 
 \equiv 0$.
In particular, each point of $\partial \Omega$ is a nodal critical point of $w$.
By \cite{Cheng}, $w$ must be identically zero.

Hence we may assume that $v$ is constant multiple of $u$. 
This implies that $\iota(u)$ is a constant multiple of $u$.
Because $\iota$ is an isometry we may conclude that $\iota(u) = \pm u$.
If $\iota(u) = - u$ then $u$ would vanish at some point on $\partial \Omega$.
Since $u|_{\partial \Omega} \equiv 1$, we conclude that $\iota(u) = u$.

A straightforward computation now shows that $R_c u$ is also symmetric with respect to $\iota$.
Because $\sin(m\theta), \cos(m\theta)$ are anti-symmetric with respect to $\iota$ for $m$ odd, we get 
\begin{equation}\label{eqn:vanishing_thm}
    \int_{\partial \Omega} \partial_\nu R_c \cdot \sin(m \cdot \theta) ds = 0 = \int_{\partial \Omega} \partial_\nu R_c \cdot \cos(m \cdot \theta) ds.
\end{equation}
Combining with \eqref{vanishing_final} we obtain that \eqref{eqn:vanishing_thm} holds for $m = 0, 1, 2, 3$.
Hence $\partial_\nu R_c|_{\partial \Omega}$ is orthogonal to $1, \cos(\theta), \sin(\theta), \cos{(2\theta)}, \sin(2\theta), \cos(3\theta), \sin(3\theta)$.
It follows from Strum-Liouville theory \cite{Ar} that $\partial_\nu R_c|_{\partial \Omega}$ has at least eight zeros.
\end{proof}

\section{Nodal sets of Rotational derivatives}
Let $\Omega$ be a simply-connected domain in the plane with analytic boundary $\partial \Omega$.
Let $v$ be a non-constant Neumann eigenfunction of $\Omega$ with eigenvalue $\mu$.
We recall that this means
\begin{equation*}
\Delta v =  \mu \cdot v, ~ ~ \textrm{and} ~~ \partial_\nu v \equiv 0 ~~ \textrm{along} ~~ \partial \Omega.
\end{equation*}
Here, as before, $\partial_\nu$ denotes the unit outward normal vector field along $\partial \Omega$.
Recall that, because $\partial \Omega$ is analytic, we have a neighborhood $N(\Omega)$ where $v$ extends analytically and satisfies $\Delta v = \mu \cdot v$ on $N(\Omega)$ \cite{TZ}.

Now we consider rotational derivatives $R_q v$ of $v$.
In particular, we would like to understand the nodal set $\mathcal{Z}(R_q v)$ for different choices of $q$.
Because $R_q$ is a Killing field we have $\Delta R_q v = \mu \cdot R_q v$ on all of $N(\Omega)$.
Hence, by \cite{Cheng}, $\mathcal{Z}(R_q v)$ is a locally finite graph.
Therefore, by considering a smaller neighborhood $N'(\Omega)$, if necessary, we assume that $\mathcal{Z}(R_q v) \cap N'(\Omega)$ is a finite graph.

For $c \in \partial \Omega$ we consider a coordinate system on the plane such that $c$ is the centre of this coordinate system and the $x$-axis is tangent to $\partial \Omega$ at $c$.
We refer to such a coordinate {\it adapted} to $c$ and $\Omega$.

\begin{lem}\label{lem:special_form}
Assume that, in an adapted coordinate in a neighborhood of $c$, $v$ has the expression   
\begin{equation*}
    v(x, y) = v_{0, 0} +  v_{0, 2} \cdot y^2 + p(x, y)
\end{equation*}
where $p$ is an analytic function of vanishing order at least three \footnote{This means that $Lp (0, 0) = 0$ for any second order differential operator $L$.}.
If $v_{0, 2} \neq 0$ and $q$ does not lie on the normal to $\partial \Omega$ at $c$ then $c$ is a degree two vertex of $\Zcal(R_q v)$.
\end{lem}
\begin{proof}
Let $q = (x_q, y_q)$ such that $q$ does not lie on the normal to $\partial \Omega$ at $c$.
Hence, $x_q \neq 0.$
Now, $R_q(x, y) = (x-x_q) \cdot v_y(x, y) - (y-y_q) \cdot v_x(x, y)$.
Using the expression for $v$ we have $$R_q v(x, y) = - 2 v_{0, 2} \cdot x_q \cdot y + f(x, y),$$
where $f$ is an analytic function of vanishing order at least two.
Because $\nabla R_q v$ does not vanish at $c$, the claim follows.
\end{proof}
Our next result is a bit more general.
Let $w$ satisfies $\Delta w = \lambda \cdot w$ on a domain $U$ for some $\lambda \neq 0$.
Let $\crit(w)$ denote the set of critical points of $w$.
For a $C^1$ curve $\gamma$ and $a \in \gamma$, let $N_\gamma(a)$ denote the normal line to $\gamma$ at $a$.
Finally, for a $C^2$ curve $\gamma$ we denote the radius of curvature and the centre of curvature of $\gamma$ at a point $a \in \gamma$ by $\rho_C(a)$ and $\xi_\gamma(a)$ respectively.
\begin{prop}
\label{prop:degree_along_normal}
Let $C \subset U$ be an analytic arc that is neither a subset of a circle nor a subset of a straight line.
Assume that $C \subset \crit(w)$.
%Let $p$ be a point on $C$ with nonzero curvature. 
Then the following holds
\begin{enumerate}
    \item if $q \notin N_C(p)$ then $p$ is a degree two vertex of $\Zcal(R_q w)$,

    \item if $q \in N_C(p)$, and, $q \neq \xi_C(p)$ then $p$ is a vertex of $\Zcal(R_q w)$ of degree four,

    \item if $q = \xi_C(p)$ then $p$ is a vertex of $\Zcal(R_q w)$ of degree at least six.
\end{enumerate}
\end{prop}
\begin{proof}
Because $w$ satisfies $\Delta w = \lambda \cdot w$ it follows that $w$ is real analytic. We now consider an expression for $w$ in a neighborhood of $q \in C$
\begin{equation*}
w(x, y) = w_{0, 0} + w_{1, 0} \cdot x + w_{0, 1} \cdot y + w_{2, 0} \cdot x^2 + w_{1, 1} \cdot xy + w_{0, 2} \cdot y^2 + O(3).
\end{equation*}
Here the coordinate system adapted to $C$ and $p$.

Since $C \subset \crit(w)$, it follows that $w(t) = w_{0, 0}$ for each $t \in C$. 
Since $C$ is not a single point, by \cite{Cheng}, it follows that $w_{0, 0} \neq 0$.
Because $p$ is a critical point of $w$, we have $w_{1, 0} = 0 = w_{0, 1}$.
Moreover, as $C \subset \crit(w)$ we have $C \subset \Zcal(\partial_x w)$. 
In particular, $w_{2, 0} = 0$. 
Using $\partial_y$ instead of $\partial_x$, in the above argument, we may further conclude that $w_{1, 1} = 0$.
Finally, since $\Delta w = \lambda \cdot w$, and, the linear term in the above expression of $w$ vanishes identically, it follows that the degree three polynomial in the expression for $w$ is harmonic.
In sum 
\begin{equation*}
w(x, y) = w_{0, 0} + w_{0, 2} \cdot y^2 + w_{3, 0} \cdot (x^3 -3 x y^2) + w_{0, 3} \cdot (y^3 -3 x^2 y) + O(4).
\end{equation*}
%
%Now we consider $R_q u$ near $p$.
%Since $C$ is not a subset of a circle, a straightforward computation shows that $\Zcal(R_q u)$ has at least two Cheng arcs crossing at $p$. In particular, the degree of $p$ as a vertex of $\Zcal(R_q u)$ is at least three.
%Since $C$ is analytic, the number of points $q$ on $C$ such that $N_C(q)$ is parallel to $N_C(p)$ is finite. We may assume, after replacing $U$ by a smaller open set, that $U$ contains no $q \neq p$ such that $N_C(q)$ is parallel to $N_C(p)$. 
%For $q \in U \cap C$, let $G(q) = (0, w(q))$ denote the point of intersection of $N_C(p)$ and $N_C(q)$.
Now consider a point $q \in N_C(p)$. In the adapted coordinates $q= (0, y_q).$ A straightforward computation shows that 
\begin{equation*}
R_q w (x, y) = 2 \cdot w_{0, 2} \cdot xy + y_q \cdot (3 \cdot w_{3, 0} \cdot (x^2 - y^2) - 6 \cdot w_{0, 3} \cdot xy) + O(3).
\end{equation*}
Since $y_q \neq 0$, if $w_{3, 0} \neq 0$ then $\Zcal(R_q w)$ would have two sub-arcs that passes through $p$ and, both of these arcs would intersect the $x$ axis transversally. 
This is impossible because $C \subset \Zcal(R_q w)$.
Hence, $w_{3, 0} = 0$ and we have the following expression
\begin{equation}\label{eqn:first_expression}
w(x, y) = w_{0, 0} + w_{0, 2} \cdot y^2 + w_{0, 3} \cdot (y^3 -3 x^2 y) + O(4).
\end{equation}
The first part of the proposition now follows from Lemma \ref{lem:special_form}.
For the rest of the claims we begin with the following.
\begin{lem}
$w_{0, 3} = 0$ if and only if $\rho_C(p) = 0$.
\end{lem}
\begin{proof}
%
%Observe that $u$ is an analytic function near $C$ and $C \subset \Zcal(\partial_x u) \cap \Zcal(\partial_y u)$.
Let $y = f_p(x)$ be the equation of $C$ near $p$ where $f_p$ is an analytic function of $x$. 
Since the coordinate system is adapted to $C$ and $p$ we have $f_p(0) = 0= f'_p(0)$.
Now, using Weierstrass preparation theorem we get 
\begin{equation}\label{eqn:second_expression}
u(x, y) =  u_{0, 0} + (y -f_p(x))^2 \cdot g_p(x, y)
\end{equation}
where $g_p(x, y)$ is an analytic function of $x, y$.

Now we compare the expressions \eqref{eqn:first_expression} and \eqref{eqn:second_expression} for $w$.
Since $f_p(0) = 0$, we get $g_p(0, 0) = w_{0, 2}$.
Since $w_{0, 0} \neq 0$ and $\Delta w = \lambda \cdot w$ we have $w_{0, 2} \neq 0$. 
Hence, $g_p(0, 0) \neq 0$.
Since $f'_p(0) = 0$, comparing the coefficient of $x^2 y$ in \eqref{eqn:first_expression} and \eqref{eqn:second_expression} we get $$2f''_p(0) \cdot w_{0, 2} = 3w_{0, 3}.$$
Since $ f''_p(0) = 1/2 \cdot \rho_C(p)$, the claim follows.
\end{proof}
\begin{remk}\label{remk:radius-of-curvature}
It follows from the proof that $\rho_C(p) = \frac{3 w_{0, 3}}{w_{0, 2}}$. Hence $\xi_C(p) = (0, \frac{w_{0, 2}}{3w_{0, 3}})$ in the chosen coordinates.
\end{remk}
First assume that $\rho_C(p) \neq 0$, and hence, $\xi_C(p)$ is finite.
In particular, $w$ has the following expression near $p$
\begin{equation*}
w(x, y) = w_{0, 0} + w_{0, 2} \cdot y^2 + w_{0, 3} \cdot (y^3 -3 x^2 y) + O(4),
\end{equation*}
where each of the three quantities $w_{0, 0}, w_{0, 2}$ and $w_{0, 3}$ is nonzero. 
Consider $q \in N_C(p)$. As our coordinates are adapted to $C$ and $p$, we have $q = (0, y_q)$. 
Therefore,
\begin{equation*}
R_q w(x, y) = (2w_{0, 2} - 6 y_q \cdot w_{0, 3}) \cdot xy + O(3).
\end{equation*}
If $q \neq \xi_C(p)$ then, by Remark \ref{remk:radius-of-curvature}, $y_q \neq 1/3 \cdot (w_{0, 2}/w_{0, 3})$.
Hence there are exactly two arcs in $\Zcal(R_q w)$ crossing each other at $p$.
This proves that $p$ is a degree four vertex of $\Zcal(R_q w)$.
This proves the second claim.

Finally, if $q = \xi_C(p)$ then $y_q = 1/3 \cdot (w_{0, 2}/w_{0, 3})$.
Therefore, the order of vanishing of $R_q w$ at $p$ is at least three.
This means that there are at least three arcs in $\Zcal(R_q w)$ crossing each other at $p$.
Hence, the degree of $p$ as a vertex of $\Zcal(R_q w)$ is at least six.
\end{proof}

\section{Proof of main theorem}
Let $\Omega$ be a centrally symmetric convex domain with analytic boundary $\partial \Omega$ and with centre $c$.
Let $u$ be a solution of \eqref{eqn:Schiffer}.
We consider $R_c u$ - the rotational derivative of $u$ with respect to the center $c$ of $\Omega$.
Our goal is to determine the points $p$ on the boundary $\partial \Omega$ such that the normal derivative $\partial_\nu R_c u(p) = 0.$ 
Because $R_c u$ vanishes along $\partial \Omega$ we obtain that $\partial_\tau R_c u|_{\partial \Omega} = 0$.
Therefore, if $\partial_\nu R_c u(p) = 0$ for some point $p \in \partial \Omega$, then such a $p$ must be a nodal critical point of $R_c u$.

Because $R_c u$ is a Laplace eigenfunction and $\partial \Omega \subset \Zcal(R_c u)$, it follows from \cite{Cheng} that as a graph, the degree of $\Zcal(R_c u)$ is at least four.
By Proposition \ref{prop:degree_along_normal} this means that $c$ lies on the normal to $\partial \Omega$ at $p$.
In particular, because $\tau(\Omega)$ - the number of points $q$ on $\partial \Omega$ such that the normal line to $\partial \Omega$ at $q$ passes through $c$ is $< 8$, we conclude that the number of points $p \in \partial \Omega$, where $\partial_\nu R_c u(p) = 0$, is at most $7$.
This contradicts Theorem \ref{thm:deng}.

%%%%%%%%%%%%%%%%%%%%%%%%%%%%%%%%%%%%%%%%%%%%%%%%%


\begin{thebibliography}{99}
%%%%%%%%%%%%%%%%%%%%%%%%%%%%%%%%%%%%%%%%%%%%%%%%%
\bibitem[Ar04]{Ar} V. I. Arnold, {\em Lectures on Partial Differential Equations.} Universitext, 2004. Springer-Verlag.

\bibitem[Av86]{Av} P. Aviles, {\em Symmetry theorems related to Pompeiu's problem.} Amer. J. Math. 108 (5) (1986) 1023-1036.

\bibitem[Ber80]{Ber} C. A. Berenstein, {\em An inverse spectral theorem and its relation to the Pompeiu problem.} J. Analyse Math. 37 (1980), 128-144.

\bibitem[BY82]{BY} C. A. Bernstein, P. Yang, {\em An overdetermined Neumann problem in the unit disc.} Adv. in Math., 4 (1982), 1-17.

\bibitem[BK82]{BK} L. Brown and J. P. Kahane, {\em A note on the Pompeiu problem for convex domains.} Math.
Ann. 259(1982), 107-110.

\bibitem[BST73]{BST} L. Brown, B. M. Schreiber and A. B. Taylor, {\em Spectral synthesis and the Pompeiu problem.} Ann. Inst. Fourier (Grenoble) 23 (1973), 125-154.

\bibitem[Cha]{Cha} I. Chavel, {\em Eigenvalues in Riemannian Geometry.
Including a chapter by Burton Randol. With an appendix by Jozef Dodziuk.}
Pure Appl. Math., 115. Academic Press, Inc., Orlando, FL, 1984. 

%\bibitem[Ca77]{Ca} L. A. Caffarelli, {\em The regularity of free boundaries in higher dimensions.} Acta Math. 139 (1977), 155-184.

%\bibitem[Cha44]{Cha} L. Chakalov, {\em Sur un problème de D. Pompeiu.} Ann. Univ. Sofia Fac. Phys. Math. 40 ( 1944), 1-44.

\bibitem[Chn76]{Cheng} S. Y. Cheng, {\em Eigenfunctions and nodal sets.} 
Comment. Math. Helv. 51 (1976), no. 1, 43-55.

\bibitem[D]{D} R. Dalmasso, {\em A new result on the Pompeiu problem.}
Trans. Amer. Math. Soc. 352 (2000), no. 6, 2723–2736.

\bibitem[De12]{De} J. Deng, {\em Some results on the Schiffer conjecture in $R^2$.} 
J. Differential Equations 253 (2012), no. 8, 2515–2526.

\bibitem[D-L15]{D-L} G. Domokos, Z. L\'{a}ngi, 
{\em On the average number of normals through points of a convex body.} (English summary)
Studia Sci. Math. Hungar.52(2015), no.4, 457–476.

\bibitem[E93]{E1} P. Ebenfelt, {\em Some results on the Pompeiu problem.}
Ann. Acad. Sci. Fenn. Ser. A I Math. 18 (1993), no. 2, 323–341.

\bibitem[E93']{E2} P. Ebenfelt, {\em Singularities of solutions to a certain Cauchy problem and an application to the Pompeiu problem.} Duke Math. J. 71 (1993) 119–142.

\bibitem[E94]{E3} P. Ebenfelt, {\em Propagation of singularities from singular and infinite points in certain complex analytic Cauchy problems and an application to the Pompeiu problem.} Duke Math. J. 73 (1994) 561–582.

\bibitem[GF91]{GF1} N. Garofalo, F. Segala, {\em New results on the Pompeiu problem.}
Trans. Amer. Math. Soc. 325 (1991), no. 1, 273–286. 

\bibitem[GF91']{GF2} N. Garofalo, F. Segala, {\em Asymptotic expansions for a class of Fourier integrals and applications to the Pompeiu problem.} J. Analyse Math. 56 (1991), 1–28. 

\bibitem[GF93]{GF3} N. Garofalo, F. Segala, {\em Another step toward the solution of the Pompeiu problem in the plane.}
Comm. Partial Differential Equations 18 (1993), no. 3-4, 491–503. 

\bibitem[GF94]{GF4} N. Garofalo, F. Segala, {\em Univalent functions and the Pompeiu problem.}
Trans. Amer. Math. Soc. 346 (1994), no. 1, 137–146.

\bibitem[KM20]{KM} B. Kawohl, L. Marcello, {\em Some results related to Schiffer's problem.} J. Anal. Math. 142 (2020), no. 2, 667–696. 

\bibitem[L07]{L} G. Liu, {\em Symmetry theorems for the overdetermined eigenvalue problems.} J. Differential Equations 233 (2007), no. 2, 585–600.

%\bibitem[Ljs59]{Lojasiewicz} S . Lojasiewicz, {\em Sur le probl\'{e}me de la division}, Studia Math. 18 (1959), 87 – 136. MR 0107168 106, 107

\bibitem[NSY20]{NSY} N. Nigam; B. Siudeja; B. Young, {\em A proof via finite elements for Schiffer's conjecture on a regular pentagon.} 
Found. Comput. Math. 20 (2020), no. 6, 1475–1504. 

\bibitem[TZ09]{TZ} Toth, John A.; Zelditch, Steve; {\em Counting nodal lines which touch the boundary of an analytic domain.} J. Differential Geom.81(2009), no.3, 649–686.

%\bibitem[Otl-Rss09]{Otal-Rosas} J.P. Otal and E. Rosas, {\em Pour toute surface hyperbolique de genre $g$, $\lambda_{2g-2} > 1/4.$} Duke Math. J. 150 (2009), no. 1, 101–115, MR 2560109, Zbl 1179.30041.

%\bibitem[Pom29]{P1} D. Pompeiu, {\em Sur certains systèmes d'équations linéaires et sur une propriété intégrale des fonctions de plusieurs variables.} C. R. Acad. Sei. Paris 188 (1929), 1.138-1.139.

%\bibitem[Pom29']{P2} -, {\em Sur une propriété intégrales des fonctions de deux variables réelles.} Bull. Sei. Acad. Royale Belgique 15 (1929), 265-269.

\bibitem[Will76]{W1} S. A. Williams, {\em A partial solution of the Pompeiu problem.} Math. Ann. 223 (1976), 183-190.

\bibitem[Will81]{W2} S. A. Williams, {\em Analyticity of the boundary for Lipschitz domains without the Pompeiu property.} Indiana Univ. Math. J. 30 (1981), 357-369.

\bibitem[Y]{Y92} S. T. Yau, {\em Seminars on Differential Geometry,} Annals of Math. Princeton Univ. Press. 1992.
  
\end{thebibliography}
\end{document}